\documentstyle[twoside,12pt]{article}

\newtheorem{th}{Theorem}[section]
\newtheorem{pro}[th]{Proposition}

\newtheorem{cor}[th]{Corollary}
\newtheorem{prb}[th]{Problem}
\newcommand{\de}{{\it des}} 
\newcommand{\mja}{{{\it maj}_A}} 
\newcommand{\dn}{{\it ndes}} 
\newcommand{\mn}{{\it nmaj}}
\newcommand{\fm}{{\it fmaj}}
\newcommand{\fd}{{\it fdes}}  
\newcommand{\inv}{{\it inv}}
\newcommand{\D}{{\it Des}}
\newcommand{\DN}{{\it NDes}}
\newcommand{\N}{{\it Neg}}
\newcommand{\n}{{\it neg}}
\begin{document}
\begin{center}
{\Large \bf Descent Numbers and Major Indices \\
for the Hyperoctahedral Group}
 \vspace{0.8cm} \\
 Ron M. Adin \footnote{
Research of all authors was supported in part by the Israel Science Foundation,
administered by the Israel Academy of Sciences and Humanities.}\\
 Department of Mathematics and Computer Science \\
 Bar-Ilan University \\
 Ramat-Gan 52900, Israel\\
 radin@math.biu.ac.il
 \vspace {0.4cm} \\
 Francesco Brenti $^1$ \\
Dipartimento di Matematica \\
Universit\'{a} di Roma ``Tor Vergata''\\
Via della Ricerca Scientifica \\
00133 Roma, Italy\\
brenti@mat.uniroma2.it
 \vspace{0.4cm} \\
Yuval Roichman $^1$ \\
 Department of Mathematics and Computer Science \\
Bar-Ilan University \\
 Ramat-Gan 52900, Israel\\
yuvalr@math.biu.ac.il 
 \vspace{0.4cm} \\
December 1, 2000
 \vspace{0.4cm} \\
{\it Dedicated to Dominique Foata}
\\
\end{center}

\newpage

\begin{abstract}
We introduce and study three new statistics on the hyperoctahedral group
$B_{n}$, and show that they give two generalizations of Carlitz's 
identity for the descent number and major index over $S_{n}$. This 
answers a question posed by Foata.
\end{abstract}

\section{Introduction}

A well known classical result due to MacMahon (see \cite{MM}) 
asserts that
the inversion number and major index of a permutation are equidistributed
over the symmetric group $S_{n}$. 
The joint distribution of major index and descent number was studied by
Carlitz~\cite{Ca}, Gessel~\cite{Ge}, and others.
Despite the fact that an increasing number of enumerative
results of this nature have been generalized to the hyperoctahedral group
$B_{n}$ (see, e.g., \cite{Bre94}, \cite{FH}, \cite{Rei93c}, \cite{Rei95}, 
\cite{Sta76}) and that several ``major index''
statistics have been introduced and studied for $B_{n}$ (see, e.g., 
\cite{CF94}, \cite{CF95a}, \cite{CF95b}, \cite{Rei93a}, \cite{Rei93b}, 
\cite{Ste94}), no generalization of MacMahon's result has been 
found until the recent paper \cite{AR}.  There a new statistic, the
{\it flag major index} (denoted here $\fm$)
was introduced; and it was shown to be equidistributed with length,
which is the natural analogue of inversion number from a Coxeter group theoretic
point of  view. No corresponding ``descent statistic'' has been found that 
allows the generalization to $B_{n}$ of the well known Carlitz identity for
descent number and major index on $S_{n}$ (see~\cite{Ca}, and also Theorem~\ref{2.3}
below),  a problem first posed by Foata  \cite{Fopc}.

\smallskip

\begin{prb} (Foata) :
Extend the (``Euler-Mahonian") bivariate distribution of descent number and
major index to the hyperoctahedral group $B_n$.
\end{prb}

\medskip

The purpose of this paper is to introduce and study three new  statistics
on $B_{n}$: the {\it negative descent}
(denoted $\dn$); the {\it negative major} 
($\mn$); and the {\it flag descent}
($\fd$). When restricted to $S_{n}$ they reduce to descent number,
major index, and twice descent number, respectively,
and they
solve the above problem. More precisely, we show that $\mn$ is
equidistributed with length on $B_{n}$, and that $\dn$ is the ``right''
corresponding descent statistic needed to extend Carlitz's result to $B_{n}$, thus 
answering Foata's question.
Finally, we prove the surprising result  that the pair of statistics
($\fd$, $\fm$) is equidistributed with  ($\dn$, $\mn$) over  $B_{n}$, thus obtaining  
a second generalization of Carlitz's identity to $B_{n}$.

\smallskip

The organization of the paper is as follows. In the next section we collect
several  definitions, notation, and results that will be used in the sequel.  In section
3 we introduce  a new ``descent set'' for elements of $B_{n}$ and,
 correspondingly,
a new ``descent number'' $\dn$  and a new ``major index'' $\mn$.
It is shown that $\mn$ is equidistributed with length over $B_{n}$ 
(Proposition \ref{3.1}), and, moreover, that the pair $(\dn,\mn)$ gives a generalization 
of Carlitz's identity (Theorem~\ref{3.2}) and thus solves Foata's problem. 
In section 4 we introduce a ``flag analogue'' of the descent number
for $B_{n}$, $\fd$, and show that it gives a second solution to Foata's problem
(Corollary \ref{4.2}). We then deduce that the pairs of statistics  $(\dn,\mn)$
and
$(\fd, \fm)$ are equidistributed over $B_{n}$ (Corollary \ref{4.4}).

\section{Notation, Definitions, and Preliminaries}

In this section we collect  some definitions,
notation and results that will be used in the rest of this paper.
We let {\bf P} $\stackrel{\rm def}{=} \{ 1,2,3, \ldots \} $,
{\bf N}$\ \stackrel{\rm def}{=}\ ${\bf P} $ \cup\ \{ 0 \} $, 
$\bf  Z$ be the ring of integers, and $\bf Q$ be the field of rational numbers;
for $a \in ${\bf N} we  let $[a] \stackrel{\rm def}{=}
\{ 1,2, \ldots , a \} $ (where $[0] \stackrel{\rm def}{=}
\emptyset $). Given $n, m \in {\bf Z}$, $n \leq m$, we let $[n,m]
\stackrel{\rm def}{=} \{ n,n+1, \ldots ,m \}$. 
For $S \subset {\bf N}$ we 
 write $S= \{ a_{1}, \ldots , a_{r} \} _{<}$ to mean that $S= \{ a_{1},
\ldots , a_{r} \}$ and $a_{1} < \ldots < a_{r} $.
The cardinality of a set $A$ will be denoted by  
$|A|$. More generally, given
 a multiset $M =  \{ 1^{a_{1}},2^{a_{2}},\ldots , r^{a_{r}} \}$
we denote by $|M|$  its cardinality, so $|M|=\sum_{i=1}^{r}a_{i}$. Given a
statement $P$ we will sometimes find it convenient to let
\[
\chi (P) \stackrel{\rm def}{=} \left\{ \begin{array}{ll}
1, &  \mbox{if $P$ is true,} \\
0, &  \mbox{if $P$ is false.}
\end{array} \right. 
\]
Given a variable $q$ and a commutative ring $R$ we denote by $R[q]$ 
(respectively,  $R[[q]]$) the ring of polynomials (respectively,
formal power series) in $q$ with coefficients in $R$. For $i \in {\bf N}$
we let, as customary, $[i]_{q} \stackrel{\rm def}{=} 1+q+q^{2}+
\ldots + q^{i-1}$ (so $[0]_{q}=0$).

 Given a sequence $\sigma =(a_{1}, \ldots ,a_{n}) \in $ {\bf
Z}$^{n}$ we say that a pair $(i,j) \in [n] \times [n]$ is an 
{\em inversion} of $\sigma $ if $i<j$ and $a_{i}>a_{j}$.
We say that $i \in [n-1]$ is a {\em descent} of $\sigma $ if
$a_{i}>a_{i+1}$. We denote by $\inv(\sigma )$ (respectively, 
$\de(\sigma )$) the number of inversions (respectively, descents) of
$\sigma $. We also let
\[ maj (\sigma ) \stackrel{\rm def}{=} \sum _{ \{ i : \; a_{i} > a_{i+1} \}
}  i \]
and call it the {\em major index} of $\sigma$.

Given a set $T$ we let $S(T)$ be the set of all bijections
$\pi : T \rightarrow T$, and $S_{n} \stackrel{\rm def}{=}S([n])$.
If $\sigma \in S_{n}$ then we write $\sigma = \sigma
_{1} \ldots \sigma _{n}$ to mean that $\sigma (i) =\sigma _{i}$,
for $i=1, \ldots ,n$. If $\sigma \in S_{n}$ then we
 will also write $\sigma $ in {\em disjoint
cycle form} (see, e.g., \cite{StaEC1}, p.~17), and will usually omit
to write the 1-cycles of $\sigma $. For example, if $\sigma =36
5492187$ then we also write $\sigma =(9,7,1,3,5)(2,6)$. Given
$\sigma , \tau \in S_{n}$ we let $\sigma  \tau \stackrel{\rm def}{=}
\sigma \circ \tau $ (composition of functions) so that, for example,
$(1,2)(2,3)=(1,2,3)$.

We denote by $B_{n}$ the group of all bijections $\pi$ of the set $[-n,n]
\setminus \{ 0 \}$ onto itself such that 
\[ \pi (-a)=-\pi (a) \]
for all $a \in [-n,n] \setminus \{ 0 \}$, with composition as the group operation.
 This group is usually known as the group of ``signed permutations''
 on $[n]$, or as the {\em hyperoctahedral group} of rank $n$. We identify $S_{n}$
 as a subgroup of $B_{n}$, and $B_{n}$ as a subgroup of $S_{2n}$, in the
 natural ways.

If $\pi \in B_{n}$ then we write $\pi = [a_{1}, \ldots ,a_{n}]$ to mean that
$\pi (i)=a_{i}$ for $i=1, \ldots , n$, and we let
\[ \inv(\pi ) \stackrel{\rm def}{=} \inv (a_{1}, \ldots , a_{n} ) ,\] 
\[ \de_{A}(\pi ) \stackrel{\rm def}{=} \de(a_{1}, \ldots , a_{n}) ,\] 
\begin{equation}
\label{2.0.3}
\mja (\pi ) \stackrel{\rm def}{=} maj (a_{1}, \ldots , a_{n} ) , 
\end{equation}
\[ \N(\pi ) \stackrel{\rm def}{=} \{ i \in [n]: \; a_{i}<0 \} , \]
and 
\[ \n (\pi ) \stackrel{\rm def}{=} |\N(\pi )| . \]
It is well known (see, e.g., \cite[Proposition 8.1.3]{BB}) that $B_{n}$
is a Coxeter group with respect to the generating set $\{ s_{0},s_{1},s_{2},
\ldots , s_{n-1} \}$, where
\[ s_{0} \stackrel{\rm def}{=} [-1,2, \ldots n ] \]
and
\[ s_{i} \stackrel{\rm def}{=} [1,2,\ldots ,i-1,i+1, i,i+2, \ldots n ] \]
for $i=1, \ldots , n-1$. This gives rise to two other natural statistics on $B_{n}$
(similarly definable for any Coxeter group), namely
\[ l(\pi ) \stackrel{\rm def}{=} min \{ r \in {\bf N}: \; \pi  =s_{i_{1}}
\ldots s_{i_{r}}  \; \; \mbox{ for some } \; \; i_{1}, \ldots ,i_{r} \in [0,n-1]
\} \]
(known as the {\em length} of $\pi$) and
\[ \de_{B}(\pi ) \stackrel{\rm def}{=} |\{ i \in [0,n-1]: l(\pi s_{i}) < 
l(\pi ) \} | . \]
There is a well known direct combinatorial way to compute these two 
statistics for $\pi \in B_{n}$ (see, e.g., \cite[Propositions 8.1.1 and 
8.1.2]{BB} or \cite[Proposition 3.1 and Corollary 3.2]{Bre94}), namely
\begin{equation}
\label{2.0.7}
 l(\pi) = inv (\pi ) - \sum _{ i \in \N(\pi )} \pi (i)  
 \end{equation}
and
\begin{equation}
\label{2.0.4}
 \de_{B}(\pi ) = | \{ i \in [0,n-1]: \; \pi  (i) > \pi  (i+1) \} | ,
 \end{equation}
where  $\pi (0) \stackrel{\rm def}{=} 0$. 
For example, if $\pi = [-3, 1,-6,2,-4,-5] \in B_{6}$ then 
$inv (\pi )=9$, $\de_{A}(\pi)=3$, 
$\mja(\pi )=11$, $\n(\pi ) =4$, $l(\pi )=27$, and $\de_{B}(\pi )=4$.

Let
\begin{equation}
\label{2.0.5}
T \stackrel{\rm def}{=} \{\pi \in B_{n}: \; \de_{A}(\pi )=0 \} .
\end{equation}
It is then well known, and easy to see, that
\begin{equation}
\label{2.0.6}
B_{n} = \biguplus _{u \in S_{n} } \{ \pi u : \; \pi \in T \} ,
\end{equation}
where $\biguplus$ denotes disjoint union. 

We will use this decomposition often in this paper.
The reader familiar with Coxeter groups will immediately recognize that 
(\ref{2.0.6})
is one case of the multiplicative decomposition of a Coxeter group into a
parabolic subgroup and its minimal coset representatives (see, e.g.,
\cite{Hum} or \cite{BB}).

\medskip

As customary, given a variable $t$ we define an operator $\delta _{t}: \;
{\bf Q}[q,t] \rightarrow {\bf Q}[q,t]$ by 
\[ \delta _{t}(P(q,t)) = \frac{P(q,qt)-P(q,t)}{qt-t} , \]
for all $P \in  {\bf Q}[q,t]$. Note that $\delta _{t}(q^{n})=0$ and
\begin{equation}
\label{2.1.0}
 \delta _{t} (t^{n}) =[n]_{q} \, t ^{n-1} 
 \end{equation}
for all $n \in {\bf N}$, and
\begin{equation}
\label{2.1.1}
\delta _{t}(A(q,t) \, B(q,t))=\delta _{t}(A(q,t))B(q,t)+A(q,qt) \delta _{t}(B(q,t))
\end{equation}
for all $A,B \in {\bf Q}[q,t]$. Also,
\[ \delta _{t} (P)(1,t) = \frac{d}{dt}(P(1,t)) \]
for all $P \in {\bf Q}[q,t]$.

For $n \in {\bf P}$ we let
\[ A_{n}(t,q) \stackrel{\rm def}{=} \sum _{\sigma \in S_{n}}t^{\de_{A}(\sigma )
} q^{\mja(\sigma )} , \]
and  $A_{0}(t,q) \stackrel{\rm def}{=} 1$. For example, $A_{3}(t,q)=1+
2tq^{2}+2tq+t^{2}q^{3}$. The following two results are due to Carlitz
\cite{Ca} and Gessel \cite{Ge}, and proofs of them can also be found in 
\cite{Gar79}.
\begin{th}
\label{2.2}
Let $n \in {\bf P}$. Then 
\[
A_{n}(t,q) =(1+tq[n-1]_{q}) \, A_{n-1}(t,q)+tq(1-t)\delta _{t}(A_{n-1}
(t,q)) . \]
\end{th}
\begin{th}
\label{2.3}
Let $n \in {\bf P}$. Then
\[
\sum _{r \geq 0} [r+1]_{q}^{n}t^{r} = \frac{A_{n}(t,q)}{{\displaystyle
\prod_{i=0}^{n}} (1-tq^{i})} . 
\]
in ${\bf Z}[q][[t]]$.
\end{th}

\section{The ``Negative" Statistics}

In this section we define and study a new ``descent set'' for the 
elements  of $B_{n}$. This gives rise, in a very natural way, to new
``major index'' and ``descent number'' statistics for $B_{n}$. We then
show that these two statistics give a generalization of Carlitz's identity
to $B_{n}$, and that one of them is equidistributed with length.

\subsection{The Negative Descent Multiset}

For $\pi \in B_{n}$ we define
\[ \D_{A}(\pi ) \stackrel{\rm def}{=}  \{ i \in [n-1]: \; \pi (i) > \pi (i+1)
\}  \]
and the {\it negative descent multiset}
\begin{equation}
\label{3.0.2}
 \DN(\pi ) \stackrel{\rm def}{=} \D_{A}(\pi ) \biguplus \{ - \pi (i) : \;
i \in \N(\pi ) \} , 
\end{equation}
where $\N(\pi)$ is the set of positions of negative entries in $\pi$, 
defined in (\ref{2.0.3}).

For example, if $\pi = [-3,1,-6,2,-4,-5] \in B_{6}$ then $\D_{A}(\pi  )
= \{ 2,4,5 \}$ and $\DN(\pi ) = \{ 2,3,4^{2},5^{2},6 \}$.
Note that if $\pi \in S_{n}$ then $\DN (\pi)$ is a set and coincides with the
usual descent set of $\pi$.
Also note that $\DN(\pi )$ can be defined rather naturally also in purely
Coxeter group theoretic terms. In fact, for $i \in [n]$ let $\eta _{i}
\in B_{n}$ be defined by
\[ \eta _{i} \stackrel{\rm def}{=}  [1,2, \ldots , i-1,-i,i+1, \ldots ,n], \]
so $\eta _{1} =s_{0}$. Then $\eta _{1}, \ldots , \eta _{n}$ are reflections
(in the Coxeter group sense; see, e.g., \cite{Hum}) of $B_{n}$ and it is clear
from (\ref{2.0.7}) that
\[
 \DN (\pi ) = 
\{ i \in [n-1]: \; l(\pi s_{i})< l(\pi ) \}
\biguplus  
\{ i \in [n]: \; l(\pi ^{-1}\eta _{i} ) < l(\pi ^{-1} ) \}  . 
\]
These considerations explain why it is natural to think of $\DN(\pi )$ as 
a ``descent set''. With this in mind, the following definitions are also 
natural. For $\pi \in B_{n}$ we let
\[ \dn (\pi ) \stackrel{\rm def}{=}  |\DN(\pi )| \]
and 
\[ \mn (\pi ) \stackrel{\rm def}{=} \sum _{i \in \DN(\pi )} i . \]
For example, if $\pi = [-3,1,-6,2,-4,-5] \in B_{6}$ then
$\dn(\pi ) = 7$ and $\mn (\pi ) =29$.

Note that  from (\ref{3.0.2}) it follows that
\begin{equation}
\label{3.0.9}
\mn(\pi )=\mja(\pi )- \sum _{i \in \N(\pi )}\pi (i) \qquad (\forall \pi\in B_n)
. 
\end{equation}
This formula is also one of the motivations behind our definition of
$\mn (\pi)$, because of the corresponding formula (\ref{2.0.7}).

\subsection{Equidistribution and Generating Function}

Our first result shows that $\mn$ and  $l$ are equidistributed in $B_{n}$.
\begin{pro}
\label{3.1}
Let $n \in {\bf P}$. Then
\[ \sum _{\pi \in B_{n}}q^{nmaj(\pi)} = \sum _{\pi \in B_{n}}q^{l(\pi)}. \]
\end{pro}
{\bf Proof.} Let $T$ be defined  by (\ref{2.0.5}). 
It is clear from our definitions that 
for all $u \in S_{n}$ and $\sigma \in T$,
$$
\mja(\sigma u)=\mja(u), \ \ \ \ inv(\sigma u)=inv (u) , 
$$
 and 
\[ \sum _{i \in \N(\sigma u)}(\sigma u)(i) = \sum _{i \in \N(\sigma )}
\sigma (i). \] Therefore, from (\ref{2.0.7}), (\ref{2.0.6}),
(\ref{3.0.9}), and the corresponding classical result for $S_n$ (see, e.g., \cite{MM}) we conclude that

\bigskip

\begin{eqnarray*}
\sum _{\pi \in B_{n}}q^{\mn (\pi )} & =  & \sum _{\sigma \in T} 
\sum _{u \in S_{n}} 
q^{\mn ( \sigma u)} \\
& = & \sum _{\sigma\in T} \sum _{u \in S_{n}} q^{\mja(\sigma u)- \sum_{i\in
\N(\sigma u)} (\sigma u)(i)} \\
& = &  \sum _{\sigma \in  T} q^{- \sum _{i \in \N(\sigma )} \sigma (i)}
\cdot \sum _{u \in S_{n}} q^{\mja(u)} \\
& = & \sum _{\sigma \in  T} q^{- \sum _{i \in \N(\sigma )} \sigma (i)}
\cdot \sum _{u \in S_{n}} q^{inv(u)} \\
& = & \sum _{\sigma \in  T} \sum _{u \in S_{n}} q^{inv( \sigma u) - 
\sum _{i \in \N (\sigma u)} (\sigma u) (i) } \\
& = & \sum _{\pi \in B_{n}} q^{l(\pi )} , 
\end{eqnarray*}
as desired. $\Box$

We now prove the main result of this section, which is also the first main result
of this work, namely that the pair of statistics $(\dn,\mn)$ gives a 
generalization of Theorem \ref{2.3} to $B_{n}$.
\begin{th}
\label{3.2}
Let $n \in {\bf P}$. Then 
\[ 
\sum _{r \geq 0} [r+1]_{q}^{n} t^{r} = \frac{
{\displaystyle \sum _{\pi \in B_{n}}} t^{ndes(\pi )} q^{nmaj(\pi )} 
}{ (1-t) {\displaystyle \prod_{i=1}^{n}} (1-t^{2}q^{2i}) }
\]
in ${\bf Z}[q][[t]]$.
\end{th}
{\bf Proof.} 
Let $T$ have the same meaning as in (\ref{2.0.5}).
 Then it is clear from
our definitions that $\de_{A}( \sigma u) =\de_{A}(u)$ for all $\sigma \in T$
and $u \in S_{n}$. Therefore we have from (\ref{2.0.6}) that
\begin{eqnarray*}
\sum _{\pi \in B_{n}}  t^{\dn(\pi )} q^{\mn (\pi )} & = & 
\sum _{\sigma \in T} \sum _{u \in S_{n}} 
t^{\de_{A}(\sigma u)+\n(\sigma u)}
q^{\mja(\sigma u)- \sum _{i \in \N (\sigma u)} (\sigma u)(i)} \\
& = & \sum _{\sigma \in T}t^{\n(\sigma )}q^{-\sum _{i \in \N(\sigma )}
\sigma (i)}  \; \cdot \; \sum _{u \in S_{n}}t^{\de_{A}(u)}q^{\mja(u)} \\
& =  & \sum _{S \subseteq [n]}t^{|S|} \; q ^{\sum _{i \in S} i} \; \cdot \;
\sum _{u \in  S_{n}} t^{\de_{A}(u)} q^{\mja(u)} \\
& = & \prod_{i=1}^{n} (1+tq^{i}) \;\cdot\; \sum _{u \in S_{n}}
t^{\de_{A}(u)} \, q^{\mja(u)} 
\end{eqnarray*}
and the result follows from  Theorem \ref{2.3}. $\Box$

\section{The Flag Statistics}

In this section we introduce yet another pair of statistics
on $B_{n}$, and show that it also gives a solution to Foata's problem.
We then derive some consequences of this result.

\subsection{The Flag Major Index}

For $i=0, 1, \ldots  ,n-1$ let $t_{i} \stackrel{\rm def}{=}s_{i}s_{i-1}
 \cdots   s_{0}$. Explicitly,
\[ t_{i} = [-i-1,1,2,\ldots ,i,i+2,\ldots ,n ] \]
for $i=0, \ldots , n-1$. It is not hard to show (see also \cite{AR})
 that
for any $\pi \in  B_{n}$ there exist unique integers $k_{0}, \ldots ,
k_{n-1}$, with $0 \leq k_{i} \leq 2i+1$ for $i=0, \ldots , n-1$, such that
\[ \pi = t_{n-1}^{k_{n-1}} \ldots t_{2}^{k_{2}}t_{1}^{k_{1}}t_{0}^{k_{0}}. \]
The {\em flag major index} of $\pi$ (see \cite{AR}) is then defined by
\[ \fm (\pi ) \stackrel{\rm def}{=} \sum _{i=0}^{n-1} k_{i} . \]
There is a simple way to compute the flag major index of a signed permutation
$\pi$. In fact, we have the following result which was 
first proved in \cite{AR}.
\begin{th}
\label{2.1}
Let $\pi \in B_{n}$. Then
\[ \fm (\pi ) =2 \, \mja(\pi ) + \n(\pi ) . \]
\end{th}
For example, if $\pi = [ -3,1,-6,2,-4,-5]$ then $\fm(\pi )=2 \cdot 11+4=26$.

\smallskip

In \cite{AR} it was shown that $\fm$ appears naturally in the Hilbert series of (diagonal action) invariant algebras; note that this property is not shared by $\mn$.

\subsection{The Flag Descent Number}

For $\pi \in B_{n}$ let
\begin{equation}
\label{4.0.0}
 \fd(\pi ) \stackrel{\rm def}{=} 2 \, \de_{A}(\pi )+\varepsilon _{1}(\pi ), 
\end{equation}
where $\de_A$ is as in Section 2 and
\begin{equation}
\label{4.0.2}
 \varepsilon _{1} (\pi ) \stackrel{\rm def}{=} \left\{
\begin{array}{ll}  
1, & \mbox{if $\pi  (1)<0$,} \\
0, & \mbox{otherwise.}
\end{array} \right. 
\end{equation}
This definition is motivated by Theorem \ref{2.1}, 
and because of this similarity
we call $\fd(\pi )$ the {\em flag descent number} of  $\pi$. For example,
if $\pi =[-3,1,-6,2,-4,-5]$ then $\fd(\pi)=2 \cdot 3 +1 =7$.
Note that from (\ref{2.0.4}), (\ref{2.0.3}) and (\ref{4.0.2}) it follows immediately that
\begin{equation}
\label{4.0.3}
\fd(\pi ) =\de_{A}(\pi )+\de_{B}(\pi ) 
\end{equation}
and also
\[  \fd(\pi ) =\de(\pi (-n), \ldots , \pi (-1),\pi  (1), \ldots , \pi  (n)) \]
for all $\pi \in B_{n}$.

Our aim is to show that the pair of statistics
($\fd,\fm$) gives a solution to Foata's problem and 
thus has the same joint distribution,
over $B_{n}$, as the pair $(\dn,\mn)$ defined in the previous section. 

\subsection{Main Theorems}

The main result of this section is

\begin{th}
\label{4.2}
Let $n \in {\bf P}$. Then
\[ \sum_{r \geq 0}[r+1]_{q}^{n}t^{r} = \frac{
{\displaystyle \sum_{\pi \in B_{n}}} t^{\fd(\pi)}q^{\fm(\pi)}}{
 (1-t){\displaystyle \prod_{i=1}^{n}} (1-t^{2}q^{2i})}  
\]
 in ${\bf Z}[q][[t]]$.
 \end{th}

This solves Foata's problem and implies, together with Theorem \ref{3.2}, that the two pairs
 of statistics
$(\fd,\fm)$ and $(\dn,\mn)$ are equidistributed over $B_{n}$.

\medskip

Theorem \ref{4.2} will be proved
in several steps.

\smallskip

Let, for convenience,
\begin{equation}
\label{4.0.4}
 S_{n}(t,q) \stackrel{\rm def}{=} \sum _{\pi \in B_{n}} t^{\fd(\pi )}
q^{\fm(\pi )} 
\end{equation}
for all $n \in {\bf P}$, and set $S_{0}(t,q) \stackrel{\rm def}{=}1$.
For example, $S_{1}(t,q)=1+tq$ and 
$S_{2}(t,q)=1+2tq+tq^{2}+t^{2}q^{2}+2t^{2}q^{3}+t^{3}q^{4}$.

Our first result gives a recursion satisfied by the polynomials $S_{n}(t,q)$.
\begin{th}
\label{4.1}
Let $n \in {\bf P}$. Then
\[
S_{n}(t,q) = (1+tq+t^{2}q^{2}[2n-2]_{q})S_{n-1}(t,q) 
+ tq(1-t)(1+tq) \delta _{t}(S_{n-1}(t,q)). \]
\end{th}
{\bf Proof.} 
The result is clear for $n \leq 2$, so fix $n \geq 3$. For $\sigma \in B_{n-1}$ 
and $i \in [n]$ let 
$$
\sigma _{i}
 \stackrel{\rm def}{=}
[\sigma (1), \ldots , \sigma(i-1),n,\sigma (i),\ldots , \sigma (n-1)]
$$ 
and
$$\sigma _{-i} \stackrel{\rm def}{=}
[\sigma (1), \ldots , \sigma(i-1),-n,\sigma (i),\ldots , \sigma (n-1)].
$$
Then clearly
\begin{equation}
\label{4.1.0}
S_{n} (t,q) = \sum _{\sigma \in B_{n-1}} \sum _{i=1}^{n}
\left[t^{\fd(\sigma _{i})} \; q^{\fm(\sigma _{i})} + t^{\fd(\sigma _{-i})}
\; q^{\fm (\sigma_{-i})} \right] . 
\end{equation}
Using (\ref{4.0.3})
and our definitions it is not hard to check (though with some patience)  that, 
for all $\sigma \in B_{n-1}$,
\[ \fd(\sigma _{1})=\fd(\sigma )+1 + \chi [\sigma (1)>0] ; \]
\[ \fd(\sigma _{-1})=\fd(\sigma ) + \chi [\sigma (1)>0] ; \]
\[ \fd(\sigma _{\pm i})=\fd(\sigma )+2  \chi [\sigma (i-1)<\sigma (i)] \]
for $i=2, \ldots , n-1$;
\[ \fd (\sigma _{n})=\fd (\sigma ) ; \]
and
\[ 	 \fd (\sigma _{-n})=\fd (\sigma )	+2. \]
Similarly, let
\[
d_{i}(\sigma ) \stackrel{\rm def}{=}
|\{ j \in \D_{A}(\sigma ): \; j \geq i \} | \qquad (i=1,\dots,n-1)
\]
be the number of descents in $\sigma$ from position $i$ on. Then
\[ maj_{A}(\sigma _{i})=
   maj_{A}(\sigma )+d_{i}(\sigma )+(i-1) \; \chi [\sigma(i-1)< \sigma (i)]+1 \]
and
\[ maj_{A}(\sigma _{-i})=
   maj_{A}(\sigma )+d_{i}(\sigma )+(i-1) \; \chi [\sigma(i-1)< \sigma (i)]  \]
for $i=1, \ldots , n-1$; 
\[ maj_{A}(\sigma _{n})=maj_{A}(\sigma ) , \]
and 
\[ maj_{A}(\sigma _{-n})=maj_{A}(\sigma )+n-1 .\]
Therefore, using Theorem \ref{2.1},
\begin{eqnarray}
\lefteqn{ \sum _{i=1}^{n} \left[t^{\fd(\sigma _{i})} \, q^{\fm(\sigma _{i})}+
t^{\fd(\sigma _{-i})} \, q^{\fm (\sigma _{-i})} \right] =}\qquad\qquad && 
\nonumber \\
&=& \lefteqn{ t^{\fd (\sigma ) +\chi [\sigma (1)>0]} \, q^{\fm (\sigma ) 
+2d_{1}(\sigma )+1}  (tq+1) +} 
\nonumber \\
&+& \lefteqn{ \sum _{\{ i \in [2,n-1]: \; \sigma (i-1)> \sigma (i) \} }
t^{\fd (\sigma )} \, q^{\fm (\sigma )+2d_{i}(\sigma )+1}(q+1) +} 
\nonumber \\
&+& \sum _{\{ i \in [2,n-1]: \; \sigma (i-1)< \sigma (i) \} }
t^{\fd (\sigma )+2} \, q^{\fm (\sigma )+2d_{i}(\sigma )+2i-1}(q+1) + 
\nonumber\\
\label{4.1.5}
&+& \lefteqn{ t^{\fd(\sigma )} \, q^{\fm(\sigma )} (1+t^{2}q^{2n-1}).}
\end{eqnarray}
From the definition of $d_{i}(\sigma )$ it is clear that 
\begin{equation}
\label{4.1.6}
\sum _{\{ i \in [2,n-1]: \; \sigma (i-1)> \sigma (i) \} } q^{2d_{i}(\sigma )} = 
\sum_{k=1}^{\de_A(\sigma)} q^{2(k-1)}
\end{equation}
for all $\sigma \in B_{n-1}$. On the other hand, let 
$$
\{ i_{1},i_{2},
\ldots , i_{a} \}_{<} \stackrel{\rm def}{=} \{ i \in [2,n-1]: \;  \sigma (i-1)
< \sigma (i) \} ,
$$ 
so that $a=n-2-\de_{A}(\sigma )$. Then 
\[
\sigma (i_{a})> \sigma(i_{a}+1)> \ldots >\sigma (n-1),
\]
and therefore 
\[ 
i_{a}+d_{i_{a}}(\sigma ) = i_{a}+(n-1-i_{a})=n-1 . 
\]
Similarly, for each $j \in  [a-1]$ one has 
\[
\sigma (i_{j})>\sigma (i_{j}+1)> \ldots > \sigma (i_{j+1}-1) < \sigma(i_{j+1}),
\] 
and therefore 
\[ 
i_{j}+d_{i_{j}}(\sigma) = i_{j}+d_{i_{j+1}}(\sigma )+(i_{j+1}-i_{j}-1) =
i_{j+1}+d_{i_{j+1}}(\sigma )-1. 
\]
This shows that
\begin{equation}
\label{4.1.7}
\sum _{\{ i \in [2,n-1]: \; \sigma (i-1)< \sigma (i) \} } q^{2(d_{i}(\sigma ) +i)} =
\sum_{k=1}^{a} q^{2(n-k)}
\end{equation}
for all $\sigma \in B_{n-1}$.

From (\ref{4.1.5}), (\ref{4.1.6}), and (\ref{4.1.7}) we conclude that
\begin{eqnarray*}
\lefteqn{ \sum _{i=1}^{n} 
\left[t^{\fd (\sigma _{i})} \, q^{\fm (\sigma _{i})}+
t^{\fd (\sigma _{-i})} \, q^{\fm (\sigma _{-i})} \right] =} \qquad\qquad && \\ 
&=& t^{\fd(\sigma)+ \chi [\sigma (1)>0]} \, q^{\fm(\sigma ) +2\de_{A}
(\sigma )+1}(1+tq) + \\
&+& t^{\fd(\sigma )}q^{\fm(\sigma )+1} 
(1+q) \sum_{k=1}^{\de_A(\sigma)} q^{2(k-1)} + \\
&+& t^{\fd(\sigma )+2}q^{\fm(\sigma )-1} 
(1+q) \sum_{k=1}^{a} q^{2(n-k)} + \\
&+& t^{\fd(\sigma )} q^{\fm(\sigma )}(1+t^{2}q^{2n-1}) = \\
&=& t^{\fd(\sigma)}q^{\fm(\sigma)}
\left\{ t^{\chi[\sigma (1)>0]}q^{2\de_{A}(\sigma)+1}(1+tq)+ \right. \\
&+& \sum_{j=1}^{2\de_A(\sigma)} q^j + t^{2} \sum_{j=2n-2a-1}^{2n-2} q^j + 
\left. (1+t^{2}q^{2n-1}) \right\} = \\
&=& t^{\fd(\sigma)}q^{\fm(\sigma)}
\left\{ t^{\chi[\sigma (1)>0]}q^{2\de_{A}(\sigma)+1}(1+tq)+ \right. \\
&+& \sum_{j=0}^{2\de_A(\sigma)} q^j + t^{2} \sum_{j=2n-2a-1}^{2n-1} \left. q^j 
 \right\}.
\end{eqnarray*}
Using the value of $a$, (\ref{4.0.0}) and some case-by-case analysis (depending on
the sign of $\sigma(1)$) we get 
\begin{eqnarray}
\lefteqn{ \sum _{i=1}^{n} 
\left[t^{\fd (\sigma _{i})} \, q^{\fm (\sigma _{i})}+
t^{\fd (\sigma _{-i})} \, q^{\fm (\sigma _{-i})} \right] =} \qquad\qquad && \\ 
&=& t^{\fd(\sigma )}q^{\fm(\sigma )}
\left\{ [\fd(\sigma )+1]_{q}+tq^{\fd(\sigma )+1}+ \right.
\nonumber\\
&+& \left. t^{2}([2n]_{q} - [\fd(\sigma )+2]_{q}) \right\} =
\nonumber\\
&=& t^{\fd(\sigma )}q^{\fm(\sigma )} 
\left\{ 1+q[\fd(\sigma )]_{q}+tq+tq(q-1)[\fd(\sigma)]_{q}+ \right.
\nonumber\\
\label{4.1.8}
&+& \left. t^{2}q^{2}([2n-2]_{q}-[\fd(\sigma )]_{q}) \right\}. 
\end{eqnarray}
Substituting (\ref{4.1.8}) into (\ref{4.1.0}) we now obtain 
\begin{eqnarray*}
S_{n}(t,q) &=& (1+tq+t^{2}q^{2}[2n-2]_{q})S_{n-1}(t,q) + \\
&+& (q+tq (q-1)-t^{2}q^{2}) \sum _{\sigma \in B_{n-1}}
[\fd(\sigma )]_{q}\,
 t^{\fd(\sigma )} q^{\fm(\sigma )},
\end{eqnarray*}
and the result follows from (\ref{4.0.4}) and (\ref{2.1.0}). $\Box$

Using the previous theorem we can now prove the following result.
\begin{th}
\label{4.3}
Let $n \in {\bf P}$. Then
\[ \sum_{\pi \in B_{n}}t^{\fd(\pi)}q^{\fm(\pi )} = \prod _{i=1}^{n} (1+tq^{i}) 
\cdot \sum _{\sigma \in S_{n}}t^{\de_{A}(\sigma )} q^{maj_{A}(\sigma )} .
\]
\end{th}
{\bf Proof.} 
We show that both sides satisfy the same recursion; equality of initial conditions is
easily checked ($n \leq 2$). Let, for convenience,
\[ 
\tilde{S}_{n}(t,q) \stackrel{\rm def}{=} 
\prod_{i=1}^{n}(1+tq^{i}) \cdot A_{n}(t,q) 
\]
be the right-hand side of the formula in Theorem \ref{4.3}.

It is easily verified that
\[ 
\delta _{t} \left( \prod_{i=1}^{n} (1+tq^{i}) \right) = q[n]_{q}
\prod_{i=2}^{n}(1+tq^{i}), 
\]
so that by (\ref{2.1.1})
\[
(1+tq)\delta_{t}(\tilde{S}_{n}(t,q)) = q[n]_{q} \tilde{S}_{n}(t,q) + 
\prod_{i=1}^{n+1} (1+tq^i) \cdot \delta_{t}(A_{n}(t,q)).
\]
By Theorem~\ref{2.2} we therefore conclude that
\begin{eqnarray*}
\tilde{S}_{n}(t,q) &=& \prod_{i=1}^{n} (1+tq^{i}) \cdot A_n(t,q) = \\
&=& \prod_{i=1}^{n} (1+tq^{i}) \cdot
\{(1+tq[n-1]_{q})A_{n-1}(t,q)+tq(1-t)\delta_{t}(A_{n-1}(t,q))\} = \\
&=& (1+tq^{n})(1+tq[n-1]_{q})\tilde{S}_{n-1}(t,q) + \\
&+& tq(1-t)\{(1+tq)\delta_{t}(\tilde{S}_{n-1}(t,q))-q[n-1]_{q}
   \tilde{S}_{n-1}(t,q)\} =\\
&=& (1+tq+t^{2}q^{2}[2n-2]_{q})\tilde{S}_{n-1}(t,q) 
   +tq(1-t)(1+tq) \delta _{t}(\tilde{S}_{n-1}(t,q)),
\end{eqnarray*}
and comparison with Theorem \ref{4.1} completes the proof. $\Box$

\medskip

Theorem \ref{4.3} has several interesting consequences.
The first one is that the descent statistic $\fd$ yields a solution to
Foata's problem (Theorem \ref{4.2}). 

\smallskip

{\bf Proof of Theorem \ref{4.2}.} This follows immediately from Theorems \ref{4.3}  and  
\ref{2.3}, given the definition of $A_{n}(t,q)$. $\Box$

This in turn implies, together with Theorem \ref{3.2}, that the two pairs of
statistics $(\fd,\fm)$ and $(\dn,\mn)$ are equidistributed in $B_{n}$.
\begin{cor}
\label{4.4}
Let $n \in {\bf P}$. Then
\[ \sum _{\pi \in B_{n}} t^{\dn(\pi )}q^{\mn(\pi )} =
\sum _{\pi \in B_{n}} t^{\fd(\pi )} q^{\fm (\pi )} . \]
\end{cor}

It would be interesting to have a direct combinatorial (i.e., bijective)
proof of this result. 

Finally, the special case $t=1$ of Corollary \ref{4.4},
together with Proposition \ref{3.1}, implies the following result, which extends
Theorem 2.2 of \cite{AR}.
\begin{cor}
\label{}
Let $n \in {\bf P}$. Then
\[ \sum _{\pi \in B_{n}} q^{l(\pi )} = \sum _{\pi \in B_{n}} q^{\mn(\pi )} 
= \sum _{\pi \in B_{n}} q^{\fm (\pi )} . \]
\end{cor}
{\bf Acknowledgments.} This 
paper grew out of stimulating discussions with Dominique Foata, Ira Gessel
and Doron Zeilberger during the conference ``Classical Combinatorics"
in honor of Foata's 65$^{th}$ birthday.
The work on it was carried out while the second author was visiting
the other authors at Bar-Ilan University.


\begin{thebibliography}{xx}


\bibitem{AR}
R. M. Adin and Y. Roichman, {\em The flag major index and group actions 
on polynomial rings}, Europ. J. Combinatorics, to appear.

\bibitem{BB} A. Bj\"{o}rner and F. Brenti, 
{\em Combinatorics of Coxeter Groups}, 
Graduate Texts in Mathematics, Springer-Verlag, 2001, to appear. 

\bibitem{Bre94}
F. Brenti, {\em $q$-Eulerian polynomials arising from Coxeter groups}, 
Europ. J. Combinatorics, {\bf 15} (1994), 417--441. 

\bibitem{Ca}
L. Carlitz, {\em A combinatorial property of $q$-Eulerian numbers},
Amer. Math. Monthly, {\bf 82} (1975), 51--54.

\bibitem{CF94}
R. J. Clarke and D. Foata, {\em Eulerian calculus. I. 
Univariable statistics}, 
Europ. J. Combinatorics, {\bf 15} (1994), 345--362.

\bibitem{CF95a}
R. J. Clarke and D. Foata, {\em Eulerian calculus. II. 
An extension of Han's fundamental transformation}, 
Europ. J. Combinatorics, {\bf 16} (1995), 221--252.

\bibitem{CF95b}
R. J. Clarke and D. Foata, {\em Eulerian calculus. III. 
The ubiquitous Cauchy formula}, 
Europ. J. Combinatorics, {\bf 16} (1995), 329--355.

\bibitem{Fopc}
D. Foata, personal communication, July 2000. 

\bibitem{FH}
D. Foata and G. N. Han, {\em Calcul basique des permutations 
sign\'ees. I. Longueur et nombre d'inversions}, 
Adv. Applied Math., {\bf 18} (1997), 489--509. 

\bibitem{Gar79}
A. Garsia, {\em On the ``maj'' and ``inv'' $q$-analogues of Eulerian polynomials}, 
Linear and Multilinear Algebra, {\bf 8}(1979/80), 21--34.

\bibitem{Ge}
I. M. Gessel, {\em Generating functions and enumeration of sequences},
Ph.D. Thesis, M.I.T., 1977.

\bibitem{Hum} 
J. E. Humphreys, {\em Reflection Groups and Coxeter Groups},
Cambridge Studies in Advanced Mathematics, no.~29,
Cambridge Univ. Press, Cambridge,  1990.

\bibitem{MM} 
P. A. MacMahon, {\em Combinatory Analysis}, Chelsea, New York, 1960. 
(Originally published in 2 vols. by the Cambridge Univ. Press, 1915--1916.)
 
\bibitem{Rei93a}
V. Reiner, {\em Signed permutation statistics}, 
Europ. J. Combinatorics, {\bf 14} (1993), 553--567.

\bibitem{Rei93b}
V. Reiner, {\em Signed permutation statistics and cycle type}, 
Europ. J. Combinatorics, {\bf 14} (1993), 569--579.

\bibitem{Rei93c}
V. Reiner, {\em Upper binomial posets and signed permutation statistics}, 
Europ. J. Combinatorics, {\bf 14} (1993), 581--588.

\bibitem{Rei95}
V. Reiner, {\em The distribution of descent and length in a Coxeter group}, 
Electron. J. Combinatorics, {\bf 2} (1995), R25.

\bibitem{Sta76}
R.\ P.\ Stanley, 
{\em Binomial posets, M\"{o}bius inversion, and permutation enumeration}, 
J.\ Combinatorial Theory Ser.\ A,
{\bf 20} (1976), 336--356.

\bibitem{StaEC1}
R.\ P.\ Stanley, {\em Enumerative Combinatorics}, vol.~1, 
Wadsworth and Brooks/Cole, Monterey, CA, 1986. 

\bibitem{Ste94}
E. Steingrimsson, {\em Permutation statistics of indexed permutations}, 
Europ. J. Combinatorics, {\bf 15} (1994), 187--205.

\end{thebibliography}
\end{document}